\documentclass[11pt]{article}
\usepackage{amsmath,amssymb,theorem}
\usepackage{diagbox}
\usepackage{multirow}
\usepackage{graphicx}
\usepackage{subfigure}
\usepackage{color}
\usepackage{epstopdf}
\usepackage[T1]{fontenc}
\usepackage{todonotes}

\newtheorem{thm}{Theorem}[section]
\newtheorem{conj}[thm]{Conjecture}
\newtheorem{lem}[thm]{Lemma}
\newtheorem{prop}[thm]{Proposition}
\theorembodyfont{\rmfamily}
\newtheorem{defn}[thm]{Definition}
\def\pf{\bigskip\noindent {\bf Proof.}~~}

\def\dfn#1{{\sl #1}}

\def\less{\backslash}

\def\pf{\bigskip\noindent {\emph{Proof.}}~~}

\def\qed{ \hfill $\square$}

\def\qed{  \bigskip\hfill\vrule height3pt width6pt depth2pt}
\newtheorem{claim}{Claim}

\title{The saturation number of  $K_{3,3}$}
\author{Shenwei Huang$^1$, Hui Lei$^2$, Yongtang Shi$^3$, and Junxue Zhang$^3$\footnote{The corresponding author.}    \\[3mm]
 $^1$College of Computer Science\\
Nankai University, Tianjin 300350, China\\
  Email: shenweihuang@nankai.edu.cn\\
    $^2$School of Statistics and Data Science, LPMC and KLMDASR\\Nankai University, Tianjin 300071, China\\
 Email: hlei@nankai.edu.cn\\
  $^3$Center for Combinatorics and LPMC\\Nankai University, Tianjin 300071, China\\
 Emails: shi@nankai.edu.cn; jxuezhang@163.com\\[2mm]
}

\date{}
\begin{document}

\maketitle
\begin{abstract}

A graph $G$ is called $F$-saturated if $G$ does not contain $F$ as a subgraph (not necessarily induced) but
the addition of any missing edge to $G$ creates a copy of $F$. The saturation number of $F$, denoted by $sat(n,F)$,
is the minimum number of edges in an $n$-vertex $F$-saturated graph.
Determining the saturation number of complete bipartite graphs is one of the most important problems in the study of  the saturation number. The value of $sat(n,K_{2,2})$ was shown to be $\lfloor\frac{3n-5}{2}\rfloor$ by Ollmann, and a shorter proof was later given by Tuza.
For $K_{2,3}$, there has been a series of study aiming to determine $sat(n,K_{2,3})$ over the years.
This was finally achieved by Chen who confirmed a conjecture of Bohman, Fonoberova, and Pikhurko that $sat(n, K_{2,3})= 2n-3$ for all $n\geq 5$.
 Pikhurko and Schmitt conjectured that $sat(n, K_{3,3})= (3+o(1))n$.
In this paper, for $n\geq 9$, we give an upper bound of $3n-9$ for $sat(n, K_{3,3})$, and prove that $3n-9$ is also a lower bound when the minimum degree of the $K_{3,3}$-saturated graphs is $2$ or $5$, where it is trivial when the minimum degree is greater than $5$.
\\

\end{abstract}

\section{Introduction}
\baselineskip 16pt

All graphs in this paper are finite and simple. Throughout the paper we use the terminology and notation of \cite{W2001}.
Given a graph $G$, we use $|G|$, $e(G)$,  $\delta(G)$, and  $\Delta(G)$  to denote the number of vertices, the number of edges,
the minimum degree and the maximum degree of $G$, respectively.
Let $\overline{G}$ denote the complement graph of $G$.
For any $v\in V(G)$,
let $d_G(v)$ and $N_G(v)$ denote the degree and neighborhood of $v$ in $G$, respectively, and let $N_G[v]=N_G(v)\cup\{v\}$.
We shall omit the subscript $G$ when the context is clear.
For
$A,B\subseteq V(G)$ with $A\cap B=\emptyset$, let $A\sim B$ denote that each vertex in $A$ is adjacent to each vertex in $B$ and $G[A, B]$  be the
subgraph with vertex set $A\cup B$ and edge set $E(G[A, B])=\{xy\in E(G): x\in A, y\in B\}$.
For $S\subseteq V(G)$, we denote by $G[S]$ the subgraph of $G$ induced by $S$. Let $n$ be a positive integer.
For positive integer $k$, we let $[k]=\{1,2,\ldots,k\}$.
We denote a path, a cycle, a star, and a complete graph with $n$ vertices by $P_n$, $C_n$, $S_n$, and $K_n$, respectively.
For $r\geq2$ and positive integers $s_1,\ldots,s_r$, let $K_{s_1,\ldots,s_r}$ denote the complete $r$-partite graph with
part sizes $s_1,\ldots,s_r$. Let $G$ and $H$ be two disjoint graphs.
 Denote by $G\cup H$  the  union of $G$ and $H$.  The {\it join} $G\vee H$  is the graph obtained from  $G\cup H$ by joining each vertex of $G$  to each vertex of $H$.

Given a family of graphs $\mathcal{F}$, a graph $G$ is \dfn{$\mathcal{F}$-saturated} if no member of $\mathcal{F}$
is a subgraph of $G$, but for any $e\in E(\overline{G})$, some member of $\mathcal{F}$ is a subgraph of $G+e$.
The \dfn{saturation number} of $\mathcal{F}$, denoted by $sat(n,\mathcal{F})$,
is the minimum number of edges in an $n$-vertex $\mathcal{F}$-saturated graph. Define $sat_\delta(n,\mathcal{F})$ to be the minimum number of edges in a graph with $n$ vertices and minimum degree $\delta$ that is $\mathcal{F}$-saturated.
If $\mathcal{F}=\{F\}$, we also write $sat(n,\{F\})$ and $sat_\delta(n,\{F\})$ as $sat(n,F)$ and $sat_\delta(n,F)$, respectively.

Saturation numbers were first studied in 1964 by Erd\H{o}s, Hajnal, and Moon \cite{EHM1964}, who proved that $sat(n,K_{k+1})=(k-1)n-{{k}\choose{2}}$.
Furthermore, they proved that equality holds
only for the graph $K_{k-1}\vee \overline{K_{n-k+1}}$.
In 1986, K\'{a}szonyi and Tuza in \cite{KT1986} determined $sat(n,F)$ for $F\in\{S_k, kK_2, P_k\}$, and they proved that
$sat(n,\mathcal{F})=O(n)$ for any family $\mathcal{F}$ of graphs.
Since then, there has been extensive research on saturation numbers for various graph families $\mathcal{F}$.

We now mention some results for complete multipartite graphs.
When all but at most one parts have size $1$, Pikhurko \cite{P1999} and  Chen,  Faudree, and Gould \cite{CFG2008} independently
determined the saturation number of complete multipartite graphs with sufficiently large order.
When there are at least two parts of size at least 2, the exact values were only known for $K_{2,2}$  and $K_{2,3}$.
The exact value for $K_{2,2}$ was first determined by Ollmann \cite{O1972}. Later on, a shorter proof was given by Tuza \cite{T1989}.
For $K_{2,3}$, there have been several papers aiming to determine $sat(n,K_{2,3})$ over the years.
This was finally achieved by Chen \cite{C2014} who confirmed a conjecture of Bohman,  Fonoberova, and Pikhurko  \cite{BFP2010}
that $sat(n, K_{2,3})= 2n-3$ for all $n\geq 5$.
For the case where the graph has $r$ parts and all parts have size 2,
Gould and Schmitt \cite{GS2007} conjectured that $sat(n, K_{2,\ldots,2})=\lceil((4r-5)n-4r^2+6r-1)/2\rceil$, and
they proved the conjecture when the minimum degree of the $K_{2,\ldots,2}$-saturated graphs is $2r-3$.
For general complete multipartite graphs $K_{s_1,\ldots,s_r}$ with $s_r\geq \cdots \geq s_1\geq1$,
Bohman, Fonoberova, and Pikhurko \cite{BFP2010}
determined the asymptotic bound on $sat(n, K_{s_1,\ldots,s_r})$ as $n\rightarrow \infty$.
\begin{thm}[\cite{BFP2010}]\label{Ks1sr}
 Let $r\geq2$ and $s_r\geq \cdots \geq s_1\geq1$. Define $p=s_1+\cdots+s_{r-1}-1$. Then, for all large $n$,
\begin{align*}
(p+\frac{s_r-1}{2})n-O(n^{3/4})\leq sat(n, K_{s_1,\ldots,s_r})\leq {p\choose 2}+p(n-p)+\big\lceil\frac{(s_r-1)(n-p)}{2}-\frac {s_r^2}{8}\big\rceil.
\end{align*}
 In particular, $sat(n, K_{s_1,\ldots,s_r})=(s_1+\ldots+s_{r-1}+0.5s_r-1.5)n + O(n^{3/4}).$
\end{thm}

We continue to study the saturation number for complete multipartite graphs.
In light of the known results, studying $sat(n,K_{3,3})$ is the natural next step.
In 2008, Pikhurko and Schmitt \cite{PS2008} conjectured that $sat(n, K_{3,3})=(3+o(1))n$.

In this paper, we give an upper bound on $sat(n, K_{3,3})$. Moreover, we  consider its lower bound. In particular, we determine the exact value of $sat(n, K_{3,3})$ for $6\le n\le 8$ and  provide a lower bound on $sat(n,K_{3,3})$ when the minimum degree of the $K_{3,3}$-saturated graphs is $2$ or $5$.
The main results  are the following theorems.
%
%
%
%
%
%
\begin{thm}\label{main1}
Let $n$ be a positive integer and $n\ge 6$. Then 
$sat(n, K_{3,3})\leq\begin{cases}
		2n,      & 6\leq n \leq 8,\\[2mm]
		3n-9,   & n \geq 9.\\
	\end{cases}$
\end{thm}

\begin{thm}\label{main23}
 \begin{description}
   \item[(i)] For $6\le n\le 8$, $sat(n,K_{3,3})=2n$.
   \item[(ii)] For $n\ge 9$, $sat_2(n,K_{3,3})= 3n-9$ and $sat_5(n,K_{3,3})\ge 3n-9$.
 \end{description}
 
\end{thm}

%
%
%
Let $G$ be a $K_{3,3}$-saturated graph with $n$ vertices and $n\geq 9$. If $\delta(G)\ge 6$, then $e(G)\ge 3n\geq 3n-9$. Hence, for $n\geq9$, to determine the exact value of  $sat(n,K_{3,3})$, we only need to consider $K_{3,3}$-saturated graphs with the minimum degree  at most $5$.


An outline of this paper is as follows.
To prove Theorem \ref{main1}, we construct 
an $n$-vertex $K_{3,3}$-saturated graph with $2n$ edges when $ 6\leq n \leq 8$ and $3n-9$ edges when $n\geq 9$ in Section \ref{upper}.
In Section \ref{689}, we first prove that $sat(n,K_{3,3})\geq 2n$ when $6\le n\le 8$ in Section \ref{6n8},  then we prove $sat_{\delta}(n,K_{3,3})\ge 3n-9$ when $\delta\in\{2,5\}$ in Section \ref{n9}.

\medskip

\section{Proof of  Theorem \ref{main1}}\label{upper}
\begin{figure}[htbp]
\begin{center}
\scalebox{0.8}[0.8]{\includegraphics{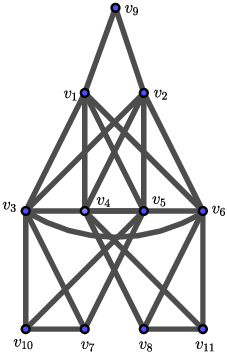}}
\caption{The graph $G_{11}$.}\label{fig1}
\end{center}
\end{figure}

In this section, for $n\geq6$, we construct an $n$-vertex $K_{3,3}$-saturated
graph $G_n$ with $2n$ edges when $ 6\leq n \leq 8$, and $3n-9$ edges when $n\geq 9$. Let $G_{11}$ be a graph as depicted in Figure \ref{fig1}. Then $G_n=G_{11}[\{v_1,\ldots,v_n\}]$ for $6\leq n\leq 11$.
\begin{prop}\label{construction1}
For $6\leq n\leq 11$, the graph $G_n$  is
$K_{3,3}$-saturated and
\begin{equation*}
e(G_n)=\begin{cases}
2n, & 6\le n\le 8,\\
3n-9, & 9\le n \le 11.\\
\end{cases}
\end{equation*}
\end{prop}
\pf
It is easy to verify that   $e(G_n)=2n$ when $6\leq n \leq 8$, and $e(G_n)=3n-9$ when $9\leq n \leq 11$. Next we show that $G_n$ contains no copy of $K_{3,3}$ for $6\leq n\leq 11$. Suppose $R$ is a copy of $K_{3,3}$ of $G_{11}$. Then $v_9\notin V(R)$ because $d_{G_{11}}(v_9)=2$. For $u\in\{v_7,v_8,v_{10},v_{11}\}$, since $d_{G_{11}}(u)=3$ and there exists $v\in N_{G_{11}}(u)$ such that $d_{G_{11}}(v)=3$ and $|N_{G_{11}}(u)\cap N_{G_{11}}(v)|=2$, we have $u\notin V(R)$. Thus $R\subseteq G_6$. Since $v_1v_2\notin E(G_6)$, $v_1$ and $v_2$ lie in the same part of $R$.  Then $R[\{v_3,v_4,v_5,v_6\}]$ contains a copy of $K_{1,3}$, a contradiction. So $G_{11}$ contains no copy of $K_{3,3}$. Note that $G_n$ ($6\leq n\leq 10$) is a subgraph of $G_{11}$.  Hence $G_n$ contains no copy of $K_{3,3}$ for any $6\leq n\leq 11$.

Let $xy$ be an edge in the complement of $G_n$.
It remains to show that the graph $G_n'$ obtained by adding $xy$ to $G_n$
has a copy of $K_{3,3}$. We consider the following cases.
\begin{description}
  \item[(a)] If $\{x,y\}\cap\{v_1,v_2\}\neq \emptyset$ or $x,y\in\{v_7,v_8,v_{10},v_{11}\}$, then the subgraph of $G_n'$
induced by $\{x,y\}\cup \{v_3,v_5\}\cup\{v_4,v_6\}$ contains a copy of $K_{3,3}$.
  \item[(b)] If $\{x,y\}\cap\{v_3,v_5\}\neq \emptyset$ or $x=v_9$, $y\in\{v_8, v_{11}\}$, then the subgraph of $G_n'$
induced by $\{x,y\}\cup \{v_1,v_2\}\cup\{v_4,v_6\}$ contains a copy of $K_{3,3}$.
  \item[(c)] If $\{x,y\}\cap\{v_4,v_6\}\neq \emptyset$ or $x=v_9$, $y\in\{v_7, v_{10}\}$, then the subgraph of $G_n'$
induced by $\{x,y\}\cup \{v_1,v_2\}\cup\{v_3,v_5\}$ contains a copy of  $K_{3,3}$.
\end{description}
For $6\leq n\leq 11$, in all cases, $G_n'$ contains a copy of $K_{3,3}$, hence $G_n$ is $K_{3,3}$-saturated.
\qed

\begin{defn}\label{Gn}
For $n\geq12$, let $H=\overline{K}_{2}\vee(C_4\cup C_{n-9}\cup K_1)$, where $V(\overline{K}_{2})=\{v_1,v_2\}$, $C_4=v_3v_4v_5v_6v_3$, $C_{n-9}=v_7v_8\ldots v_{n-3}v_7$, $V(K_1)=\{v_{n-2}\}$. Let $G_n$ be the graph obtained from $H$ by adding new vertices $\{v_{n-1},v_n\}$ and new edges $\{v_{n-1}v_3,v_{n-1}v_5,v_nv_4,v_nv_6\}$.
\end{defn}

%
%
%

\begin{prop}\label{construction}
For $n\geq12$, the graph $G_n$ defined in Definition~\ref{Gn} is
$K_{3,3}$-saturated and has $3n-9$ edges.
\end{prop}
\pf Clearly, $e(G)=2(n-4)+(n-5)+4=3n-9$. Firstly, We show that $G_n$ has no subgraph isomorphic to $K_{3,3}$. Suppose $R$ is a copy of $K_{3,3}$ of $G_n$. From the structure of  $G_n$, we see that $d(v_{n-1})=d(v_n)=2$ and hence $v_{n-1}, v_n\notin V(R)$. Thus $R\subseteq H$. Since each vertex of $C_4\cup C_{n-9}\cup K_1$ has at most two neighbors in $C_4\cup C_{n-9}\cup K_1$,  $v_1,v_2\in V(R)$ and they lie in different parts of $R$. This contradicts $v_1v_2\notin E(G_n)$. So $G_n$ contains no copy of $K_{3,3}$.

Let $xy$ be an edge in the complement of $G_n$.
It remains to show that the graph $G''$ obtained by adding $xy$ to $G_n$
has a copy of $K_{3,3}$. We consider the following cases.
\begin{description}
  \item[(a)] If $x,y\in\{v_1,v_2,v_{n-1},v_{n}\}$, then the subgraph of $G''$
induced by $\{x,y\}\cup\{v_3,v_5\}\cup\{v_4,v_6\}$ contains a copy of $K_{3,3}$.
  \item[(b)] If $x=v_{n-1}$, $y\in\{v_4,v_6,v_7,\ldots,v_{n-2}\}$ or $x=v_4$, $y=v_6$ or $x\in\{v_4,v_6\}$, $y\in\{v_7,\ldots,v_{n-2}\}$, then the subgraph of $G''$
induced by $\{x,v_1,v_2\}\cup\{y,v_3,v_5\}$ contains a copy of $K_{3,3}$.
  \item[(c)] If $x=v_{n}$, $y\in\{v_3,v_5,v_7,\ldots,v_{n-2}\}$ or $x=v_3$, $y=v_5$ or $x\in\{v_3,v_5\}$, $y\in\{v_7,\ldots,v_{n-2}\}$, then the subgraph of $G''$
induced by $\{x,v_1,v_2\}\cup\{y,v_4,v_6\}$ contains a copy of $K_{3,3}$.
 \item[(d)]If  $x,y\in \{v_7,\ldots,v_{n-2}\}$ and $x\neq v_{n-2}$, let $N(x)\cap\{v_7,\ldots,v_{n-3}\}=\{x',x''\}$,  then the subgraph of $G''$
induced by $\{x,v_1,v_2\}\cup\{y,x',x''\}$ contains a copy of $K_{3,3}$.
\end{description}
In all cases, $G''$ contains a copy of $K_{3,3}$. Hence $G_n$ is $K_{3,3}$-saturated.
\qed

By Proposition \ref{construction1} and Proposition \ref{construction}, we complete the proof of Theorem \ref{main1}.\bigskip

\section{Proof of Theorem \ref{main23}}\label{689}

In the rest of the paper, we consider the lower bound on $sat(n, K_{3,3})$. Let $G=(V,E)$ be a $K_{3,3}$-saturated graph. We firstly choose a vertex $a$ such that $d(a)=\delta(G)$ and $e(G[N(a)])$ is as small as possible. 
We partition $V$ into four parts $V_1$, $V_2$, $V_3$ and $V_4$,
where $V_1=N[a]$, $V_2=\{x \in V \backslash V_1: |N(x) \cap N(a)| \geq 2 $\},
$V_3=\{y \in V \backslash (V_1 \cup V_2): |N(y) \cap N(a)|=1$\} and $V_4=V\backslash (V_1 \cup V_2 \cup V_3 )$.
Let $N_G(a)=\{a_1,a_2,\ldots,a_{d(a)}\}$.
For $i_1,i_2,\ldots, i_s\in [d(a)]$, let $V_{i_1i_2\ldots i_s}=\{x \in V_2: N(x) \cap V_1=\{a_{i_1},a_{i_2},\ldots, a_{i_s}\}\}$.

In the following, we will first describe some useful properties of the $K_{3,3}$-saturated graph $G$.

\begin{prop}\label{k22}
 The following statements hold.
\begin{description}
\item[(i)] For any $x,y\in V$, if  $x y \notin E$, then there are $\{x_1, x_2\}\subseteq N(x) $ and $\{y_1, y_2\}\subseteq N(y)$ such that $\{x_1,x_2\}\sim\{y_1, y_2\}$. (We usually say there is a copy of $K_{2,2}$ between $N(x)$ and $N(y)$.)
\item[(ii)] For any $x \in V \setminus V_1$ , we have $|N(x) \cap N(a_i) \cap N(a_j)|\leq2$ for any $i,j\in[d(a)]$ with $i\neq j$, and there exist $i,j\in[d(a)]$ with $i\neq j$ such that $|N(x) \cap N(a_i) \cap N(a_j)|=2$.
\item[(iii)] For any $x \in V_3$, we have $|N(x) \cap V_2| \geq 1$.   For any $x \in V_4$, we have $|N(x) \cap V_2| \geq 2$.

\item[(iv)] When $G[V_1 \backslash \{a\}]$ contains no copy of $K_{1,2}$, we have $|N(x) \cap V_2|\geq2$ for any $x \in V \setminus V_1$, and  $|V_2|\geq3$.
When $G[V_1\setminus\{a\}]$ contains no copy of $K_{2,2}$, we have $|N(x)\cap V_2|\ge 1$ for any $x\in V_2$, and $|V_2|\ge 2$.
\end{description}
\end{prop}
\pf Suppose  $x y \notin E$. Then there is a copy of $K_{3,3}$ in $G+xy$, and (i) follows. For any $x \in V \setminus V_1$,
if there is a vertex $x \in V \setminus V_1$ such that $|N(x) \cap N(a_i) \cap N(a_j)|\geq3$ for some $i,j\in[d(a)]$ with $i\neq j$,
then we would obtain a copy of $K_{3,3}$ of $G$, a contradiction. So $|N(x) \cap N(a_i) \cap N(a_j)|\leq2$ for any  $x \in V \setminus V_1$ and $i,j\in[d(a)]$ with $i\neq j$. Since $ax\notin E$ for any $x\in V\less V_1$, there exist $i,j\in[d(a)]$ such that $|N(x) \cap N(a_i) \cap N(a_j)|=2$ by (i). This proves (ii).
Let $x \in V \setminus V_1$ and $i,j\in[d(a)]$ with $i\neq j$ such that $|N(x) \cap N(a_i) \cap N(a_j)|=2$, we say $\{u,v\}=N(x) \cap N(a_i) \cap N(a_j)$. Then $u,v\in (V_1\cup V_2)\setminus\{a\}$. If $x \in V_3$, then we have $|N(x) \cap V_2| \geq 1$ by the definition of $V_3$.  If $x \in V_4$, then we have $|N(x) \cap V_2| \geq 2$ by the definition of $V_4$.  This proves (iii).
Suppose $G[V_1 \backslash \{a\}]$ contains no copy of $K_{1,2}$. Then $u,v\in  V_2$. Hence we have $|N(x) \cap V_2|\geq2$ for any $x \in V \setminus V_1$,  and $|V_2|\geq3$.
Suppose $G[V_1 \backslash \{a\}]$ contains no copy of $K_{2,2}$. Then $\{u,v\}\cap  V_2\ne \emptyset$. Hence we have $|N(x) \cap V_2|\geq1$ for each $x \in V_2$,  and $|V_2|\geq2$.
 This proves (iv).
\qed

Proposition \ref{k22}(i) implies $\delta(G)\ge 2$ for each $K_{3,3}$-saturated graph $G$. Thus we consider $\delta(G)\ge 2$.



\subsection{Proof of Theorem \ref{main23}(i)}\label{6n8}
By Theorem \ref{main1}, to prove $sat(n, K_{3,3})= 2n$ for $6\leq n\leq 8$, it suffices to prove $sat(n, K_{3,3}) \geq 2n$.
We consider the minimum degree of $G$. If $\delta(G) \geq 4$,
then we have $e(G) \geq 2n$. So we assume that $2\leq \delta(G)\leq 3$. For $i\in \{2,3,4\}$ and $x\in V_i$, we define $f(x)=|N(x) \cap (V_1 \cup \dots \cup V_{i-1} )|+0.5|N(x) \cap V_i|-2$.
Let  $ s_i=\sum_{x \in V_i} f(x)$, where $i \in \{2,3,4\}$.

We first observe that one can relate the number of edges to $s_2$, $s_3$ and $s_4$ in the following way:
\begin{align}\label{2ef}
e(G)=\notag&~e(G[V_1])+e(G[V_2])+e(G[V_1,V_2])+e(G[V_3])+e(G[V_1,V_3])+e(G[V_2,V_3])+e(G[V_4])\\\notag&+e(G[V_4,V_2\cup V_3])\\\notag
=&~e(G[V_1])+2(|V_2|+|V_3|+|V_4|)+s_2+s_3+s_4\\
=&~e(G[V_1])+2(n-|V_1|)+s_2+s_3+s_4.
\end{align}

\begin{lem}\label{n678}
For $6\le n\le 8$,
\begin{description}
  \item[(i)] if $\delta(G)=2$, then $s_2+s_3+s_4\geq |V_2|+|V_3|$.
  \item[(ii)] if  $\delta(G)=3$, then $s_2+s_3+s_4 \geq |V_2|+|V_3|+|V_4|$ when $e(G[V_1 \backslash \{a\}])\leq 1$ and $s_2+s_3+s_4 \geq 0.5(|V_2|+|V_3|+|V_4|)$ when $e(G[V_1 \backslash \{a\}])\geq2$.
\end{description}
\end{lem}
\pf Suppose that $\delta(G)=2$.  Then $G[V_1 \backslash \{a\}]$ contains no $K_{1,2}$. Thus $f(x) \geq 1$ for each $x\in V_2\cup V_3$ and $f(x) \geq 0$ for each $x\in V_4$ by Proposition \ref{k22} (iii). So $s_2+s_3+s_4 \geq |V_2|+|V_3|$. Suppose that  $\delta(G)=3$.  If $e(G[V_1 \backslash \{a\}])\leq 1$, then  $|V_4|\leq1$ because $n\leq8$ and $|V_2|\ge 3$ by Proposition \ref{k22}(iv). Thus $f(x) \geq 1$ for each $x\in V\setminus V_1$ by Proposition \ref{k22} (iii). So $s_2+s_3+s_4 \geq |V_2|+|V_3|+|V_4|$. If $e(G[V_1 \backslash \{a\}])\geq2$, then we have $|N(x) \cap V_2|\geq1$ for each $x \in V_2\cup V_3$ and $|N(x) \cap V_2|\geq2$ for each $x \in V_4$ by Proposition \ref{k22} (iii). Thus for $x \in V_2$, $f(x) \geq 0.5$; for $y \in V_3$, $f(y) \geq 0.5$ or $f(y)=0$ and there exists a vertex $z \in V_4$ such that $f(z)=1$; for $z \in V_4$, $f(z) \geq 0.5$. Proposition \ref{k22}(iv) implies  $|V_2|\ge 2$ and so
  $|V_3\cup V_4|\leq 2$,  we have $s_2+s_3+s_4 \geq 0.5(|V_2|+|V_3|+|V_4|)$.
\qed


Suppose that $\delta(G)=2$. If $a_1a_2 \in E$, then $e(G)\geq2n+|V_2|+|V_3|-3$ by Lemma \ref{n678}(i) and (\ref{2ef}). By Proposition \ref{k22}(iii), we have $|V_2| \geq 3$. So $e(G) \geq 2n$.  If $a_1a_2 \notin E(G)$, then $e(G)\geq2n+|V_2|+|V_3|-4$ by Lemma \ref{n678}(i) and (\ref{2ef}). Proposition \ref{k22}(i) implies that there is a copy of $K_{2,2}$ between $N(a_1)$ and $N(a_2)$, we have $|V_2\cup V_3| \geq 4$. So $e(G) \geq 2n$.\medskip

Suppose that $\delta(G)=3$. If $n=6$, then $|V_2|=2$, $|V_3|=|V_4|=0$ and $e(V_1)=6$ by Proposition \ref{k22}(i). Otherwise, $a_ia_j\notin E$ where $i,j\in [3]$ with $i\ne j$, Proposition \ref{k22}(i) implies that there is a copy of $K_{2,2}$ between $N(a_i)$ and $N(a_j)$, which contradicts the fact that $|V_2\cup V_3|=2$.
 Let $V_2=\{x_1, x_2\}$. Proposition \ref{k22}(iv) implies $x_1x_2\in E$.  If $x_1a_i\notin E$ for some $i\in [3]$, then $x_2\in V_{123}$ by Proposition \ref{k22} (i).
Thus $e(G)\ge 12=2n$.

If $n=7$  and $e(G[V_1])\le 4$, then $G[V_1\setminus\{a\}]$ contains no copy of $K_{1,2}$. Proposition \ref{k22}(iv) implies
$|V_2|=3$
and $e(G[V_2])=3$. Since $a_ia_j\notin E$ for some $i,j\in[3]$, Proposition \ref{k22}(i)  implies there is a copy of $K_{2,2}$ between $N(a_i)$ and $N(a_j)$. Since $|V_2\cup V_3|=|V_2|= 3$,  $e(G[V_1])\ge 4$.
We see $|V_{123}|\le 1$, else $G$ contains a copy of $K_{3,3}$. There exists a vertex $x$ such that $|N(x)\cap V_1|=2$ and  $xa_k\notin E$ for some $k\in [3]$. Proposition \ref{k22}(i) implies that there is a copy of $K_{2,2}$ between $N(x)$ and $N(a_k)$, say $\{x_1, x_2\}\sim\{a_{k1},a_{k2}\}$.
When $\{a_{k1},a_{k2}\}\subseteq V_2$, then $\{x_1,x_2\}\subseteq V_1$ and $\{a_{k1},a_{k2}\}\subseteq V_{123}$, a contradiction.
When $\{a_{k1},a_{k2}\}\cap V_1\ne \emptyset$, since $e(G[V_1])\le 4$,  $|\{a_{k1},a_{k2}\}\cap\{a_1,a_2,a_3\}|\le 1$.
If $a_{k1}\in \{a_1,a_2,a_3\}$, then $a_{k2}\in V_2$. By $|V_2|=3$, $\{x_{k1},x_{k2}\}\cap V_1\ne \emptyset$, which contradicts $e(G[V_1])\le 4$.
If $a\in \{a_{k1},a_{k2}\}$, say $a_{k1}=a$, then $\{x_{1},x_2\}\subseteq V_1$, $a_{k2}\in V_2$ and $a_{k2}\in V_{123}$. Then $e(G)=e(G[V_1])+e(G[V_2])+e(G[V_1,V_2])\ge 4+3+7=14=2n$.

If $n=7$ and $e(G[V_1])=6$, by Lemma \ref{n678}(ii), then $e(G)\ge 2n-0.5$, that is $e(G)\ge 2n$.
Suppose $n=7$ and $e(G[V_1])=5$. Let $E(G[V_1\setminus\{a\}])=\{a_1a_2,a_1a_3\}$.
 If  $|V_2|=2$, then let $V_2=\{x_1,x_2\}$. Applying Proposition \ref{k22}(i) to $ax_1\notin E$ ($ax_2\notin E$), we have the $K_{2,2}$ between $N(a)$ and $N(x_1)$ ($N(x_2)$) is $\{a_2,a_3\}\sim\{a_1, x_2\}(\{a_1, x_1\})$.
 Then $\{x_1,x_2\}\subseteq V_{123}$, and so $\{a_1,a_2,a_3\}\sim\{a, x_1, x_2\}$ is a copy of $K_{3,3}$ of $G$, a contradiction.
 If $|V_2|\ge 3$, then $|V_2|=3$ by $n=7$.  Let $V_2=\{x_1,x_2, x_3\}$. Note that $f(x_i)\ge 0.5$ for each $i\in [3]$. If there exists a vertex $x_i\in V_{123}$ or there are two vertices $x_i,x_j\in V_2$ such that $f(x_i)\geq1$ and $f(x_j)\geq1$,  then $e(G)\ge 2n-0.5$ by (\ref{2ef}),  and so $e(G)\ge 2n$. Thus we may assume $V_{123}=\emptyset$ and there is at most one vertex $x_i\in V_2$ such that $f(x_i)\geq1$.
Since there is a copy of $K_{2,2}$ between $N(x)$ and $N(a)$ for each $x\in V_2$, there is some vertex $x_i\in V_2$ with $f(x_i)=1$, say $x_1$. Then $x_1\in V_{23}$ and $\{x_2,x_3\}\subseteq V_{1i}$ for some $i\in \{2,3\}$, say $i=2$. Then $N(a_3)=\{a,a_1,x_1\}$, but $e(G[N(a_3)])\le 1$, which contradicts the minimality of $e(G[N(a)])$.  So $e(G)\ge 2n$.


If $n=8$, then $e(G) \geq 2n$ when $e(G[V_1 \backslash \{a\}])=1$ or $3$ by Lemma \ref{n678}(ii). Suppose $n=8$ and $e(G[V_1 \backslash \{a\}])=0$, then $ e(G)=2n+s_2+s_3+s_4-5$. So we need to show $s_2+s_3+s_4 \geq 4.5$.
 If $|V_{123}|\geq1$, then $f(x) \geq 2$ for each $x\in V_{123}$. So $s_2+s_3+s_4 \geq |V_2|+|V_3|+|V_4|+1 \geq 5$ by the proof of Lemma \ref{n678}(ii). Now we  consider $|V_{123}|=0$.
 Since $a_1a_2\notin E$, Proposition \ref{k22}(i) implies that there  is  a copy of $K_{2,2}$ between $N(a_1)$ and $N(a_2)$, say $\{x_1,x_2\}\sim \{x_3,x_4\}$. Then $\{x_1,x_2,x_3,x_4\}\subseteq V_2\cup V_3$. Since $n=8$,  $|V_2\cup V_3|=4$. If $x_1\in V_3$, then we can not find  a copy of $K_{2,2}$ between $N(a_2)$ and $N(a_3)$ because $|(N(a_2)\cup N(a_3))\cap (V_2\cup V_3)|\leq 3$, a contradiction. By symmetry, we have $\{x_1,x_2,x_3,x_4\}\subseteq V_2$.
 If there exists $i\in[4]$ such that $|N(x_i)\cap V_2|\ge 3$, then $e(G)= e(G[V_1])+e(G[V_2])+e(G[V_1,V_2])\ge3+5+8=16=2n$. If $|N(x_i)\cap V_2|=2$ for each $i\in [4]$, then $E(G[V_2])=\{x_ix_j|i\in\{1,2\}, j\in\{3,4\}\}$. Since $x_1x_2\notin E$, Proposition \ref{k22}(i) implies that there  is  a copy of $K_{2,2}$ between $N(x_1)$ and $N(x_2)$. Note that $N(x_1)\cup N(x_2)\subseteq\{a_1,a_2,a_3,x_3,x_4\}$,  $e(G[\{a_1,a_2,a_3\}])=0$ and $x_3x_4\notin E$. So the $K_{2,2}$ between $N(x_1)$ and $N(x_2)$ must be $\{a_1,a_2\}\sim \{x_3,x_4\}$. Then $d(a_3)=1$, this contradicts $\delta(G)\geq2$.
Suppose $n=8$ and $e(G[V_1 \backslash \{a\}])=2$. Then $ e(G)=2n+s_2+s_3+s_4-3$. So we need to show $s_2+s_3+s_4 \geq 2.5$. If $f(x)\geq1$ for some $x\in  V_2$, then $s_2+s_3+s_4 \geq 2.5$ by the proof of Lemma \ref{n678}(ii). If $f(x)=0.5$ for some $x\in  V_2$, then $f(x')\geq 1$ where $\{x'\}=N(x)\cap V_2$. So $s_2+s_3+s_4 \geq 2.5$.

This completes the  proof of the lower bound on $sat(n,K_{3,3})$ for $6\leq n\leq 8$.
\qed

\subsection{Proof of Theorem \ref{main23}(ii)}\label{n9}
Note that for $n\geq9$, the minimum degree of the $K_{3,3}$-saturation graph we constructed in Section \ref{upper} with $3n-9$ edges is $2$.
Thus $sat_2(n,K_{3,3})\le 3n-9$ for $n\ge 9$. Hence, to prove $sat_2(n,K_{3,3})=3n-9$, it  suffices to prove $sat_2(n,K_{3,3})\ge 3n-9$ for $n\ge 9$. In this section, we give the lower bound of $3n-9$ for $sat_{\delta}(n, K_{3,3})$ for $\delta\in\{2,5\}$ and $n\geq9$. We first consider the case where the minimum degree of $G$ is $2$.
\subsubsection{ $\delta(G)=2$}
We prove $sat_2(n, K_{3,3})\ge 3n-9$ for $n\geq 9$ in this part.
According to the partition of $V$, we define  $h(x)=|N(x) \cap (V_1 \cup \ldots\cup V_{i-1} )|+0.5|N(x) \cap V_i|-3$ for each $ x \in V_i $  and $q_i=\sum_{x \in V_i} h(x)$ where $i\in \{2,3,4\}$. For each $x\in V$, we say that the $h$-value of $x$ is $k$ if $h(x)=k$.
	\begin{align}\label{eG2}
	e(G)=&~e(G[V_1])+e(G[V_2])+e(G[V_1,V_2])+e(G[V_3])+e(G[V_1,V_3])+e(G[V_2,V_3])+e(G[V_4])\notag\\&+e(G[V_4,V_2\cup V_3])\notag\\
	=&~e(G[V_1])+3(|V_2|+|V_3|+|V_4|)+q_2+q_3+q_4\notag\\
	=&~e(G[V_1])+3(n-|V_1|)+q_2+q_3+q_4.
	\end{align}

 By (\ref{eG2}), we have $e(G) \geq 3n-7+q_2+q_3+q_4$.
Therefore, it suffices to prove
\begin{align}\label{-2.5}
q_2+q_3+q_4 \geq -2.5.
\end{align}
By  Proposition \ref{k22}(iii), we have $|N(x)\cap V_2|=2$ for each $x\in V\setminus V_1$. So $h(z) \geq 0$ for each $z \in V_2\cup V_3 $
and  $ h(z)\geq -1$ for each $z \in V_4$. Thus,
$q_2 \geq 0$ and $q_3 \geq 0$. Therefore, to prove (\ref{-2.5}), it suffices to show   $q_4 \geq -2.5$.

Let $V_4^{-}=\{z\in V_4:h(x)<0\}=\{z_1,z_2,\ldots,z_{|V_4^{-}|}\}$ and $n^-_4(x)=|N(x)\cap V^-_4|$ for each $x\in V$.
By Proposition \ref{k22}(iii), each vertex $z \in V_4^{-}$ has exactly two neighbors in $V_2$,
so we let $N(z_i)\cap V_2=\{x_{i1},x_{i2}\}$.
Note that if $h(z_i)=-1$, then $N(z_i)=\{x_{i1},x_{i2}\}$ and so $z_i$ has no neighbor in $V_4^{-}$,
and if $h(z_i)=-0.5$, then $d(z_i)=3$ and $z_i$ has one neighbor in $V_4$, saying $N_4(z_i)=\{c_i\}$.

 For each $z_i, z_j\in V_4^{-}$ with $z_iz_j\notin E$, there is a $K_{2,2}$ between $N(z_i)$ and $N(z_j)$ by Proposition \ref{k22}(i), we define four  different types of $K_{2,2}$ as follows.  \medskip

\noindent {\bf Type 1} : $\{x_{i1}, x_{i2}\}\sim\{x_{j1}, x_{j2}\}$;\medskip

\noindent {\bf Type 2} : $\{x_{i1}, x_{i2}\}\sim\{x_{jt}, c_j\}$, where $t\in\{1,2\}$;\medskip

\noindent {\bf Type 3} : $\{x_{is}, c_i\}\sim\{x_{j1}, x_{j2}\}$, where $s\in\{1,2\}$;\medskip

\noindent {\bf Type 4} : $\{x_{is}, c_i\}\sim\{x_{jt}, c_j\}$, where $s,t\in\{1,2\}$.\medskip

If there are three vertices in $V_4$ with an $h$-value of $-1$, then there are six distinct vertices $x_1,x_2,\ldots, x_6 \in V_2$
such that $\{x_1, x_2\} \sim \{x_3, x_4\}$, $\{x_3, x_4\} \sim \{x_5, x_6\}$ and $\{x_1, x_2\} \sim \{x_5, x_6\}$.
Thus $G$ contains a copy of $K_{3,3}$ as $\{a_1, a_2, x_1\} \sim \{x_3, x_4, x_5\}$,
a contradiction. So there are at most two vertices in $V_4$ with an $h$-value of $-1$. Thus $q_4 \geq -2.5$ when $|V_4^{-}|\leq 3$.  In the following,  we assume that $|V_4^{-}|\geq 4$.

\begin{claim}\label{three -1}
There is at most one vertex in $V_4^{-}$ with an $h$-value of $-1$.	
\end{claim}
\pf Suppose that, by contradiction,  there are exactly two vertices with an $h$-value  $-1$, say $z_1$ and $z_2$.
Then $z_1 z_2 \notin E$ and the $K_{2,2}$ between $N(z_1)$ and  $N(z_2)$ is Type 1.
Since $|V_4^{-}|\geq 4$, there exists a vertex, say $z_3$, such that $d(z_3)=3$ and $z_1z_3\notin E$, $z_2z_3\notin E$.
Applying Proposition \ref{k22}(i) to $z_1z_3\notin E$ and $z_2z_3\notin E$, we obtain that there exists $b \in N(z_3)$ such that $b \sim \{x_{11}, x_{12}, x_{21}, x_{22}\}$. Then $G$ contains a copy of $K_{3,3}$ as $\{a_1, a_2, b\} \sim \{x_{11}, x_{12}, x_{21}\}$, a contradiction. Hence there is at most one vertex  in $V_4^{-}$ with an $h$-value of $-1$.
 \qed

 By Claim \ref{three -1}, if $|V_4^{-}| \leq 4$, then $q_4\geq -2.5$. So we assume that $|V_4^{-}| \geq 5$ in the following.
 \begin{claim}\label{one -1}
 If there exists a vertex in $ V_4^{-}$ with an $h$-value of $-1$, then $q_2+q_3+q_4 \geq -2.5$.	
 \end{claim}

\pf  Without loss of generality, we assume that $h(z_1)=-1$. For each $z_i\in V_4^{-}\setminus\{z_1\}$, since $z_1z_i\notin E$,
$\{x_{11},x_{12}\}\not\subseteq N(z_i)$.
We first prove that there is at most one vertex $z_i\in V_4^{-}$ such that $\{x_{11}, x_{12}\}\sim \{x_{i1}, x_{i2}\}$.
Suppose not. Then there exist two vertices, say $z_2$ and $z_3$, such that $\{x_{11}, x_{12}\} \sim \{x_{t1}, x_{t2}\}$ for each $t\in \{2,3\}$.
Since $|N(x) \cap V_2| =2$ for each $x \in V \setminus V_1$, $\{x_{21}, x_{22}\} =\{x_{31}, x_{32}\}$.
Note that $z_2z_3\in E$ for otherwise the non-edge $z_2z_3$ contradicts Proposition \ref{k22}(i).
Since $|V_4^{-}| \geq 5$, there exists a vertex, say $z_4$, such that $z_4 z_p\notin E$ for each $p \in [3]$.
By applying Proposition \ref{k22}(i) to $z_1z_4$, we have $\{x_{4i},c_4\}\sim \{x_{11}, x_{12}\}$ for some $i\in \{1,2\}$ and thus $x_{4i} \in \{x_{21}, x_{22}\}$.
Since $c_2=z_3$, there is no $K_{2,2}$ between $N(z_2)$ and $N(z_4)$, contradicting  Proposition \ref{k22}(i).
This proves the statement.

Now for $i\in \{3,4,\ldots,|V_4^{-}|\}$,  since $z_1z_i\notin E$, we have  $\{x_{11},x_{12}\}\sim \{x_{ij}, c_i\}$, where $j\in[2]$.
Applying Proposition \ref{k22}(i) to $c_iz_1\notin E$,
we have $|N(c_i)\cap V_4|\ge 3$ and thus $h(c_i)\ge 0.5$.
 Now we show that $c_i\neq c_j$ for $i,j\in \{3,4,\ldots,|V_4^{-}|\}$ with $i\neq j$.
Since $c_i\notin V_4^{-}$, we have $z_iz_j\notin E$. By Proposition \ref{k22} (i), there is a $K_{2,2}$ between
$N(z_i)$ and $N(z_j)$. By considering the $K_{2,2}$ between $N(z_k)$ and $N(z_1)$ for $k\in \{i,j\}$, we see $N(c_k)\cap V_2=\{x_{11},x_{12}\}$.
It follows that the $K_{2,2}$ between $N(z_i)$ and $N(z_j)$ must be Type 4.
So $c_i\neq c_j$. Now we have
\[
	q_4\ge h(z_1)+h(z_2)+\sum_{i=3}^{|V_4^{-}|}(h(z_i)+h(c_i))\geq -1.5.
\]
This completes the proof. \qed

By Claim \ref{one -1},  we assume $h(z) =-0.5$ for each vertex $z \in V_4^{-}$.
 If $|V_4^{-}| \leq 5$, then $q_4 \geq -2.5$. So we assume  $|V_4^{-}| \geq 6$ in the following.

 \begin{claim}\label{type1}
If $h(z) =-0.5$ for each vertex $z \in V_4^{-}$ and there exist two non-adjacent vertices in $ V_4^{-}$ satisfying the $K_{2,2}$ between  their neighborhood is Type 1, then $q_2+q_3+q_4 \geq -2.5$.	
 \end{claim}
\pf Suppose  $z_1z_2\notin E$ and the $K_{2,2}$ between $N(z_1)$ and $N(z_2)$ is Type 1. Let $U=\{z \in V_4^{-}\setminus\{z_1,z_2\}$ with $zz_1, zz_2\notin E\}$.  Since $|V_4^{-}| \geq 6$, we have  $|U|\ge 2$.
Let $z_i\in U$. By applying Proposition \ref{k22} (i) to $z_i z_1 \notin E$, there is a copy of $K_{2,2}$ between $N(z_1)$ and $N(z_i)$.  Note that $|N(v)\cap V_2|=2$ for each $v\in V\setminus V_1$.  If the $K_{2,2}$ is Type 1 or Type 3, then $\{x_{i1},x_{i2}\}=\{x_{21},x_{22}\}$. If the $K_{2,2}$ is Type 2, then $N(c_3)\cap V_2=\{x_{11},x_{12}\}$ and $x_{is}\in\{x_{21},x_{22}\}$ for some $s\in[2]$. In each case, we cannot find a $K_{2,2}$ between $N(z_2)$ and $N(z_i)$. So the $K_{2,2}$ between $N(z_1)$ and $N(z_i)$ is Type 4. Similarly, the $K_{2,2}$ between $N(z_2)$ and $N(z_i)$ is Type 4. So we have $x_{is}\in\{x_{11},x_{12}\}$ and $x_{it}\in\{x_{21},x_{22}\}$, where $\{s,t\}=[2]$, and $c_1,c_2,c_i\notin V_4^{-}$.
%
%
Hence for each $z_i, z_j \in U$, the $K_{2,2}$ between $N(z_i)$ and $N(z_j)$ is
Type 4. So $c_i \neq c_j$.
This means that for each $z\in U$, its unique neighbor $c\in V_4$ has at least 3 neighbors in $V_4\setminus V^-_4$,  so $h(z)+h(c)\ge 0$. And for any $z_i, z_j \in U$, $c_i \neq c_j$, so $q_4\ge -2$.
\qed

By Claim \ref{type1}, we suppose there are no two vertices $z_i,z_j\in V_4^-$ with $z_iz_j\notin E$ such that the $K_{2,2}$ between $N(z_i)$ and $N(z_j)$ is Type 1.
Suppose that $c \in V_4^{-}$ for each $z \in V_4^{-}$. Let $z_i,z_j \in V_4^{-}$ with $z_iz_j\notin E$. By Proposition \ref{k22}(i), Claim \ref{type1}, and $c_i, c_j \in V_4^{-}$, we may assume the $K_{2,2}$ between $N(z_i)$ and $N(z_j)$ is Type $2$. Then there is no copy of $K_{2,2}$ between $N(z_i)$ and $N(c_j)$, a contradiction. So we choose $z \in V_4^{-}$ with $c \notin V_4^{-}$ as $z_1$.
Let $A_0= \emptyset$. Let $A_\ell=\{z|z\in V_4^{-} \setminus (A_0 \cup \ldots \cup A_{\ell-1}) \text{ and the $K_{2,2}$ between $N(z_1)$ and $N(z)$ is Type $\ell$ } \}$  and  $C_\ell=\{c_i:z_i\in A_\ell\}$ for $\ell \in [4]$. By Claim \ref{type1}, we have $A_1=C_1=\emptyset$. Thus $|A_2|+|A_3|+|A_4|=|V_4^{-}|-1$.
Let  $C=\{c_1\}\cup C_2\cup C_3 \cup C_4$ and $C'=\{c_1\}\cup C_2\cup C_4$.
Note that $C_j$ and $C_k$ may intersect when $j \neq k$ and $j,k\in[4]$.



For any $z\in A_2$,  we have $c \notin V_4^{-}$ for otherwise  there is no copy of $K_{2, 2}$ between $N(z_1)$ and $N(c)$.
Thus for each  $z_i,z_j\in A_2$,  we have $z_iz_j\notin E$. Since $z_i,z_j\in A_2$, we have $N(c_i)\cap V_2=N(c_j)\cap V_2=\{x_{11},x_{12}\}$ and there exist $s,t\in[2]$ such that $x_{is}\notin\{x_{11},x_{12}\}$ and $x_{jt}\notin\{x_{11},x_{12}\}$. If the $K_{2,2}$ between $N(z_i)$ and $N(z_j)$ is
Type 2 or Type 3, then $N(c_j)\cap V_2=\{x_{i1},x_{i2}\}$ or $N(c_i)\cap V_2=\{x_{j1},x_{j2}\}$, a contradiction.  So  the $K_{2,2}$ between $N(z_i)$ and $N(z_j)$ is
Type 4. This implies that $C_2$ is a clique.
%
%
For each  two vertices  $z_i,z_j\in A_3$, we have $N(z_i)\cap V_2=N(z_j)\cap V_2$ since $|N(c_1)\cap V_2|=2$.
If  $z_iz_j\notin E$,  then the $K_{2,2}$ between $N(z_i)$ and $N(z_j)$ is Type $4$. If $z_iz_j\in E$, then $c_ic_j\in E$.  This implies that $C_3$ is a clique. Thus if $|A_3| \geq 3$, then for each $z\in A_3$,  we have $c \notin V_4^{-}$.

Let $|C_2|=p$, $|C_3\setminus C_2|=q$ and $|C_4\setminus(C_3\cup C_2)|=r$. Note that $|C|\le p+q+r+1$ and  the equation $|C|= p+q+r+1$ implies that $c_1\notin C_2\cup C_3$.
Note that \begin{align}\label{q4=-2.5}
q_4\ge&\sum_{v\in C\setminus V_4^-}h(v)+\sum_{v\in V_4^-}h(v)=\sum_{v\in C\setminus V_4^-}h(v)-0.5|V_4^-|.
\end{align}
To prove $q_4\ge -2.5$, it suffices to prove $\sum_{v\in C\setminus V_4^-}h(v)\ge0.5|V_4^-|-2.5$ by (\ref{q4=-2.5}).
Recall that $C_2$ and $C_3$ are two cliques of $G$,   $(C_2\cup C_4)\cap V_4^-=\emptyset$ and $C_3\cap V_4^-=\emptyset$ if $|A_3|\geq3$.

\noindent {\bf Case 1:} $|C_3|=|A_3|\geq 3$.

In this case, we have $(C_2\cup C_3\cup C_4)\cap V_4^-=\emptyset$. Thus $\sum_{v\in C\setminus V_4^-}h(v)=\sum_{v\in C}h(v)$.

If $C_2\cap C_3\ne \emptyset$, then
\begin{align*}
\sum_{v\in C}h(v)\ge& 2|C|+e(G[C])+0.5e(G[C,V_4^-])-3|C|\\
\ge&2|C|+{{p}\choose{2}}+{{q}\choose{2}}+q+r+0.5|V_4^-|-3|C|\\
=&{{p}\choose{2}}+{{q}\choose{2}}+q+r+0.5|V_4^-|-|C|\\
\ge &max\{0,p-1\}+max\{0,q-1\}+p+r+0.5|V_4^-|-(p+q+r+1)\\
\ge & 0.5|V_4^-|-2.
\end{align*}

If $C_2\cap C_3=\emptyset$, then $q\geq3$ and
\begin{align*}
\sum_{v\in C}h(v)\ge&\notag  2|C|+e(G[C])+0.5e(G[C,V_4^-])-3|C|\\\notag
\ge &2|C|+{{p}\choose{2}}+{{q}\choose{2}}+r+0.5|V_4^-|-3|C|\\
=&{{p}\choose{2}}+{{q}\choose{2}}+r+0.5|V_4^-|-(p+q+r+1)\\
\ge& p-1+q+r+0.5|V_4^-|-(p+q+r+1)\\
=& 0.5|V_4^-|-2.
\end{align*}
\noindent {\bf Case 2:} $|A_3|\le 2$ and $|A_2|=p\geq 3$.
\begin{align*}
\sum_{v\in C\setminus V_4^-}h(v)\ge&\sum_{v\in C'}h(v)\ge 2|C'|+e(G[C'])+0.5e(G[C',V_4^-])-3|C'|\\
\ge &2|C'|+{{p}\choose{2}}+|C_4\setminus C_2|+0.5(|V_4^-|-2)-3|C'|\\
\ge&{{p}\choose{2}}+|C_4\setminus C_2|+0.5|V_4^-|-1-(p+|C_4\setminus C_2|+1)\\
\ge& {{p}\choose{2}}-p+0.5|V_4^-|-2\\
\ge& 0.5|V_4^-|-2.
\end{align*}
\noindent {\bf Case 3:} $|A_2|\le 2$ and $|A_3|\le 2$.

 Note that $(\{c_1\}\cup C_4)\cap (\{z_1\}\cup A_4)=\emptyset$. We have
 \begin{align}\label{A4-1}
  \sum_{v\in \{c_1\}\cup C_4}h(v)\ge &2(|C_4|+1)+e(G[\{c_1\}\cup C_4])+0.5e(G[\{c_1\}\cup C_4,\{z_1\}\cup A_4])-3(|C_4|+1)\notag\\
  \ge &2(|C_4|+1)+ |C_4|+0.5(|A_4|+1)-3(|C_4|+1)=0.5(|A_4|-1).
 \end{align}
Then \begin{align*}
q_4\ge & \sum_{v\in \{c_1\}\cup C_4}h(v)+\sum_{v\in V_4^-}h(v)\\
\ge &0.5(|A_4|-1)-0.5(|A_2|+|A_3|+|A_4|+1)=-0.5(|A_2|+|A_3|)-1.
\end{align*}
Observe that $q_4\ge -2.5$ when $|A_2|+|A_3|\le 3$. Thus we just need to consider the case $|A_2|=|A_3|=2$.

Note that $C'\cap V_4^-=\emptyset$. Suppose $C_2\cap (\{c_1\}\cup C_4)\ne \emptyset$. Then $G[C']$ is a connected graph, and so $e(G[C'])\ge |C'|-1$.
We see  \begin{align*}
\sum_{v\in C\setminus V^-_4}h(v)\ge\sum_{v\in C'}h(v)\ge& 2|C'|+e(G[C'])+0.5e(G[C',V_4^-\setminus A_3])-3|C'|\\
\ge & e(G[C'])-|C'|+0.5(|V_4^-|-2)\\
\ge & |C'|-1-|C'|+0.5|V_4^-|-1\\
\ge & 0.5|V_4^-|-2.
\end{align*}
Suppose $C_2\cap (\{c_1\}\cup C_4)=\emptyset$. Let $C_2=\{c_2, c_3\}$.
If  $h(c_2)>0$ or $h(c_3)>0$, by (\ref{A4-1}), then
\begin{align*}
q_4\ge&\sum_{v\in \{c_1\}\cup C_4}h(v)+\sum_{v\in C_2}h(v)+\sum_{v\in V_4^-}h(v)\\
\ge& 0.5(|A_4|-1)+0.5-0.5(1+4+|A_4|)=-2.5.
\end{align*}

 If  $h(c_2)=h(c_3)=0$, then  $N(c_2)=\{x_{11}, x_{12}, c_3, z_2\}$ and $N(c_3)=\{x_{11}, x_{12}, c_2, z_3\}$.
Since $z_1c_2\notin E$, the $K_{2,2}$ between $N(z_1)$ and $N(c_2)$ must be Type 4, which contradicts $c_3c_1\notin E$.

In a conclusion, $q_4\ge -2.5$ and so $e(G)\ge 3n-9$. This completes the proof of the lower bound on $sat_2(n,K_{3,3})$ for $\geq 9$.\qed


\subsubsection{$\delta(G)=5$}
We prove $sat_5(n, K_{3,3})\ge 3n-9$ for $n\geq 9$ in this part. Since $\delta(G)=5$, we have $e(G)\ge 2.5n$. Then $e(G)\ge 3n-9$ when $n\le 19$.   Thus we assume $n\ge 20$ in the following.

 We define a new function $g$ as follows.
 \begin{itemize}
   \item For $x\in  V_2$, $g(x)=|N(x)\cap V_1|+0.5|N(x)\cap (V_2\cup V_3)|+0.25|N(x)\cap V_4|-3$.
   \item For $x\in V_3$, $g(x)=|N(x)\cap V_1|+0.5|N(x)\cap (V_2\cup V_3\cup V_4)|-3$.
   \item For $x\in V_4$, $g(x)=0.75|N(x)\cap V_2|+0.5|N(x)\cap (V_3\cup V_4)|-3$.
 \end{itemize}
Observe that
\begin{align}\label{eg}
e(G)=\notag&~e(G[V_1])+e(G[V_2])+e(G[V_1,V_2])+e(G[V_3])+e(G[V_1,V_3])+e(G[V_2,V_3])+e(G[V_4])\\&+e(G[V_4,V_2\cup V_3])\notag\\
=&~e(G[V_1])+3(|V_2|+|V_3|+|V_4|)+\sum_{x\in V\setminus V_1}g(x)\notag\\
=&~e(G[V_1])+3(n-|V_1|)+\sum_{x\in V\setminus V_1}g(x).
\end{align}
Note that $\delta(G)=5$. Then $g(x)\ge -0.25$ for each  $x\in V_2$ because $|N(x)\cap V_1|\ge 2$; $g(x)\ge 0$ for each $x\in V_3$ because $|N(x)\cap V_1|=1$; $g(x)\ge 0$ for each $x\in V_4$ because $|N(x)\cap V_2|\ge 2$. If there exists a vertex $x_0\in V_2$ such that $g(x_0)<0$, then $g(x_0)=-0.25$, $d(x_0)=5$, $N(x_0)\cap (V_2\cup V_3)=\emptyset$, $|N(x_0)\cap V_1|=2$ and $|N(x_0)\cap V_4|=3$. We may assume that  $N(x_0)=\{a_i,a_j,z_1,z_2,z_3\}$, where $i,j \in [5]$, $i\ne j$ and $\{z_1,z_2,z_3\}\subseteq V_4$. Since $ax_0\notin E(G)$,  Proposition \ref{k22}(ii) implies that there is a copy of $K_{2,2}$ in $G[V_1\setminus\{a\}]$. Let $s=1$ if $a_ia_j\in E$ and $s=0$ if $a_ia_j\notin E$. Thus $e(G[V_1\setminus\{a\}])\ge 4+s$. But $e(G[N(x_0)])\le 3+s$ because $N(z_i)\cap V_1=\emptyset$ for each $i\in [3]$, which contradicts the minimality of $e(G[N(a)])$. Hence, $g(x)\ge 0$ for each $x\in V\setminus V_1$ and so $\sum_{x\in V\setminus V_1}g(x)\ge 0$. When $e(G[V_1])\ge 9$, by (\ref{eg}), we have $e(G)\ge 3n-9$ .
 Thus we next consider $e(G[V_1])\leq8$. Note that $|N(x)\cap V_2|\geq 1$ for each $x\in V_2$ when $e(G[V_1])\leq8$. The following discussion is split into three cases below.

\medskip
\noindent {\bf Case 1:}  $e(G[V_1])=8$.\medskip

If  $\sum_{x\in V\setminus V_1}g(x)>0$, then $e(G)=3n-10+\sum_{x\in V\setminus V_1}g(x)>3n-10$ by (\ref{eg}) and so $e(G)\ge 3n-9$ because $e(G)$ is an integer. Next we  prove  $\sum_{x\in V\setminus V_1}g(x)>0$. If there exists a vertex $x\in V_2$ with $|N(x)\cap V_1|\ge 3$, then $g(x)>0$ and so $\sum_{x\in V\setminus V_1}g(x)>0$.
So we may assume that $|N(x)\cap V_1|=2$ for each $x\in V_2$.
Since $e(G[V_1\setminus\{a\}])=3$, there is a vertex $a_i$ such that $N(a_i)\cap N(a)=\emptyset$ or $N(a_i)\cap N(a)=\{a_j\}$ with $N(a_j)\cap N(a)=\{a_i\}$, where $i,j\in [5]$ and $i\neq j$. We denote such a vertex by $a_1$. There is a vertex $a_k$ such that $a_1a_k\notin E$ for $k\in [5]$ and $k\ne 1$.
Since $a_1a_k\notin E$,  by Proposition \ref{k22}(i),  $N(a_1)\cap(V_2\cup V_3)\ne \emptyset$.  Let  $x\in N(a_1)\cap(V_2\cup V_3)$ and $x_1\in N(x)\cap V_2$.
If $x\in V_3$, then $|N(x_1)\cap (V_2\cup V_3)|\ge 2$. If $x\in V_2$, by the choice of $a_1$, then we have $|N(x_1)\cap V_2|\ge 2$, else there is no $K_{2,2}$ between $N(x_1)$ and $N(a)$.
So $g(x_1)\ge 0.25$, which implies $\sum_{x\in V\setminus V_1}g(x)>0$. Hence $e(G)\ge 3n-9$.

\medskip
\noindent {\bf Case 2:} $e(G[V_1])=7$ and there is a copy of $K_{1,2}$ in $G[V_1\setminus\{a\}]$.\medskip

 We may assume that $E(G[V_1\setminus\{a\}])=\{a_1a_2, a_1a_3\}$. If $\sum_{x\in V\setminus V_1}g(x)>1$, by (\ref{eg}), then $$e(G)=e(G[V_1])+3(n-|V_1|)+\sum_{x\in V\setminus V_1}g(x)>7+3(n-6)+1=3n-10.$$ Since $e(G)$ is an integer, $e(G)\ge 3n-9$. Thus we just need to prove $\sum_{x\in V\setminus V_1}g(x)>1$.
Let $V_2^1=\{x\in V_2: |N(x)\cap    V_2|=1 \}$ and $V_2^2=\{x\in V_2: |N(x)\cap  V_2|\ge 2 \}$.
Let $x\in V_2^1$ and $xx_1\in E(G[V_2])$. Applying  Proposition \ref{k22}(i) to $ax\notin E(G)$, we have $x\in N(a_1)$ and $x_1\in N(a_2)\cap N(a_3)$. If $x_1\in V_2^1$,  then $x_1\in N(a_1)$ and $x\in N(a_2)\cap N(a_3)$ by $x_1a\notin E(G)$.  Thus $\{a_1, a_2, a_3\}\subseteq (N(x)\cap V_1)\cap (N(x_1)\cap V_1$). There is a copy of $K_{3,3}$ in $G$,  that is $\{a,x,x_1\}\sim\{a_1, a_2, a_3\}$, a contradiction.
  This implies that
 $e(G[V_2^1])=0$, $V_2^2\ne \emptyset$ and $|V_2|\ge 3$.
Since $a_4a_5\notin E$, there is a copy of $K_{2,2}$ between $N(a_4)$ and $N(a_5)$, say $\{x_{41},x_{42}\}\sim \{x_{51}, x_{52}\}$. Notices that  $N(a_4)\cap V_1=N(a_5)\cap V_1=\{a\}$. Thus $\{x_{41},x_{42},x_{51}, x_{52}\}\subseteq V_2\cup V_3$. For each $y\in \{x_{41},x_{42},x_{51}, x_{52}\}\cap V_3$, by Proposition \ref{k22}(i), then $|N(y)\cap V_2|\ge 2$.
By the definition of $g$-function,  for each $x\in V_2$, we have
\begin{align*}
g(x)=&|N(x)\cap V_1|+0.25|N(x)\cap (V_2\cup V_3\cup V_4)|+0.25|N(x)\cap (V_2\cup V_3)|-3\\
=& |N(x)\cap V_1|+0.25|N(x)\cap (V_2\cup V_3\cup V_4)|+0.25|N(x)\cap V_2|-3+0.25|N(x)\cap V_3|.
\end{align*}
If $x\in V_2^1$, then
\begin{align*}
g(x)\geq  2+0.25\times3+0.25\times1-3+0.25|N(x)\cap V_3|
= 0.25|N(x)\cap V_3|.
\end{align*}
If $x\in V_2^2$, then
\begin{align*}
g(x)\geq  2+0.25\times3+0.25\times2-3+0.25|N(x)\cap V_3|
= 0.25+0.25|N(x)\cap V_3|.
\end{align*}
If $|N(x)\cap V_1|\geq 3$, then
\begin{align*}
g(x)\geq  3+0.25\times2+0.25\times1-3+0.25|N(x)\cap V_3|
= 0.75+0.25|N(x)\cap V_3|.
\end{align*}
Suppose $|\{x_{41},x_{42},x_{51}, x_{52}\}\cap V_3|\ge 2$.
Then  $e(G[V_2,V_3])\ge2|\{x_{41},x_{42},x_{51}, x_{52}\}\cap V_3| \ge 4$. Note that $V_2^2\ne \emptyset$. Thus  $$\sum_{x\in V\setminus V_1}g(x)\ge \sum_{x\in V_2}g(x)\ge 0.25+\sum_{x\in V_2}0.25|N(x)\cap  V_3|=0.25+0.25e(G[V_2,V_3])\ge 1.25.$$
 Suppose $|\{x_{41},x_{42},x_{51}, x_{52}\}\cap V_3|= 1$, say  $x_{41}\in V_3$. Then $\{x_{42}, x_{51}, x_{52}\}\subseteq V_2$ and $x_{42}\in V_2^2$. We see $\{x_{51},x_{52}\}\subseteq V_2^2$ or $x_{42}\in N(a_2)\cap N(a_3)$, that is $|N(x_{42})\cap V_1|\ge 3$.
Thus $$\sum_{x\in V\setminus V_1}g(x)\ge \sum_{x\in V_2}g(x)\ge 0.75+\sum_{x\in V_2}0.25|N(x)\cap  V_3|=0.75+0.25e(G[V_2,V_3])\ge 1.25.$$
 It remains to consider the case $ \{x_{41},x_{42},x_{51}, x_{52}\} \subseteq V_2$, that is $ \{x_{41},x_{42},x_{51}, x_{52}\} \subseteq V_2^2$.
 If $V_3\ne \emptyset$, then $e(G[V_2,V_3])\ge 1$ and
 $$\sum_{x\in V\setminus V_1}g(x)\ge \sum_{x\in V_2}g(x)\ge 0.25|V_2^2|+\sum_{x\in V_2}0.25|N(x)\cap  V_3|
 \ge 1+0.25e(G[V_2,V_3])\ge
  1.25.$$
If $|N(x)\cap V_1|\ge 3$ for some $x\in V_2$, then $$\sum_{x\in V\setminus V_1}g(x)\ge \sum_{x\in V_2}g(x)\ge 0.75+0.25(|V_2^2|-1)+\sum_{x\in V_2}0.25|N(x)\cap  V_3|
 \ge1.5.$$
Next we assume that $|N(x)\cap V_1|=2$ for each $x\in V_2$ and $|V_3|=0$. Note that for each $x\in V_2^1$, let $xx_1\in E(G[V_2])$, we have $x_1\in N(a_2)\cap N(a_3)$. Thus $x_1\notin \{x_{41},x_{42},x_{51}, x_{52}\}$.
If $|V_2|\ge 5$, then $V_2^2\setminus\{x_{41},x_{42},x_{51}, x_{52}\}\ne \emptyset$. Thus $|V_2^2|\ge 5$ and $\sum_{x\in V\setminus V_1}g(x)\ge 1.25$. If $|V_2|\le 4$, that is $V_2=\{x_{41},x_{42},x_{51}, x_{52}\}$, then we have $|V_4|\ge n-|V_2|-|V_3|-6=n-10$ because $|V_3|=0$. Note that $n\geq 20$. Thus
\begin{align*}
e(G)=&e(G[V_1])+e(G[V_2])+e(G[V_1\cup V_4,V_2])+e(G[V_4])\\
\ge & 7+4+8+2|V_4|+\frac{3|V_4|}{2}>3n-9
\end{align*}






\medskip
\noindent {\bf Case 3:} $e(G[V_1])=7$ and there is no copy of $K_{1,2}$ in $G[V_1\setminus\{a\}]$ or $5\leq e(G[V_1])\leq6$.\medskip

In this case, we define a new function $g'$ as follows.
\begin{itemize}
  \item For $x\in  V_2$, $g'(x)=|N(x)\cap V_1|+0.5|N(x)\cap V_2|-3$.
  \item For $x\in V_3\cup V_4$, $g'(x)=|N(x)\cap (V_1\cup V_2)|+0.5|N(x)\cap ( V_3\cup V_4)|-3$.
\end{itemize}
We see
\begin{align}\label{egprime}
e(G)=&\notag~e(G[V_1])+e(G[V_2])+e(G[V_1,V_2])+e(G[V_3])+e(G[V_1,V_3])+e(G[V_2,V_3])+e(G[V_4])\\\notag
&+e(G[V_4,V_2\cup V_3])\\\notag
=&~e(G[V_1])+3(|V_2|+|V_3|+|V_4|)+\sum_{x\in V\setminus V_1}g'(x)\\
=&~e(G[V_1])+3(n-|V_1|)+\sum_{x\in V\setminus V_1}g'(x).
\end{align}

For each $x\in V_2$, by Proposition \ref{k22}(iv), $|N(x)\cap V_2|\ge 2$. Thus $g(x)\ge0.25$ because $d(x)\ge 5$. It follows that $\sum_{x\in V\setminus V_1}g(x)\ge 0.25|V_2|$. It suffices to consider the following two subcases. \medskip

\noindent {\bf Subcase 3.1:} $|V_2|\ge 13$ or $|V_3\cup V_4|\geq 7$\medskip

Suppose $|V_2|\ge 13$. Then $$e(G)=e(G[V_1])+3(n-|V_1|)+\sum_{x\in V\setminus V_1}g(x)\ge 5+3n-18+0.25|V_2|\ge 3n-9.75$$ and so $e(G)\ge 3n-9$ because $e(G)$ is an integer.

Suppose $|V_3\cup V_4|\geq 7$.
By Proposition \ref{k22}(iv), $|N(x)\cap V_2|\ge 2$ for each $x\in V\setminus V_1$. Thus $g'(x)\ge 0$ for each $x\in V_2$, $g'(x)\ge 1$ for each $x\in V_3$, and $g'(x)\ge 0.5$ for each $x\in V_4$.
 It follows that $$e(G)=e(G[V_1])+3(n-|V_1|)+\sum_{x\in V\setminus V_1}g'(x)\ge 5+3n-18+0.5|V_3\cup V_4|\ge 3n-9.5.$$ Since $e(G)$ is an integer, $e(G)\ge 3n-9$.\medskip

\noindent {\bf Subcase 3.2:} $|V_2|\le 12$ or $|V_3\cup V_4|\leq 6$\medskip

Since $n\geq12$, $|V_3\cup V_4|\ge 2$. We first prove the following claim.
\begin{claim}\label{v34}
If there is no copy of $K_{1,2}$ in $G[V_1\setminus\{a\}]$ and $|V_3\cup V_4|\ge 2$, then 	$\sum_{x\in V_3\cup V_4}g'(x)\ge 2$. In particular, if $|V_3|\ge 1$ or $|N(z)\cap V_2|\ge 3$ for some $z\in V_4$, then $\sum_{x\in V_3\cup V_4}g'(x)\ge 3$.
\end{claim}
\pf By the definition of $g'$-function and $\delta(G)= 5$, we have  for each $x\in V_3$, $g'(x)\ge 1$ and for each $x\in V_4$, $g'(x)\ge 0.5$. When $|V_3\cup V_4|\ge 4$, $\sum_{x\in V_3\cup V_4}g'(x)\ge 2$.
When $2\le |V_3\cup V_4|\le 3$,  for each $x\in V_3\cup V_4$, we have  $|N(x)\cap (V_1\cup V_2)|=5-(|V_3\cup V_4|-1)$. Thus $$g'(x)\ge (6-|V_3\cup V_4|)+0.5(|V_3\cup V_4|-1)-3=2.5-0.5(|V_3\cup V_4|$$   and $$\sum_{x\in V_3\cup V_4}g'(x)\ge (2.5-0.5(|V_3\cup V_4|))|V_3\cup V_4|\ge 3.$$
Next we assume that  $|V_3|\ge 1$ or $|N(z)\cap V_2|\ge 3$ for some $z\in V_4$. To prove $\sum_{x\in V_3\cup V_4}g'(x)\ge 3$, it suffices to  consider the case $|V_3\cup V_4|\ge 4$ by the above discussion.
If $|V_3\cup V_4|\ge 5$ or $|V_3|\ge 2$, then $\sum_{x\in V_3\cup V_4}g'(x)\ge 3$.
Suppose $|V_3\cup V_4|= 4$ and $|V_3|\leq 1$. Let $V_3\cup V_4=\{y_1,y_2,y_3,y_4\}$ and $\{y_1,y_2,y_3\}\subseteq V_4$.
Let  $y_4\in V_3$ or $|N(y_4)\cap V_2|\ge 3$ when $y_4\in V_4$. If $g'(y_i)\geq 1$ for some $i\in [3]$, then $\sum_{x\in V_3\cup V_4}g'(x)\ge 3$. So we assume $g'(y_i)=0.5$ for each $i\in [3]$, then we have $|N(y_i)\cap (V_3\cup V_4)|=3$. Thus  $G[\{y_1,y_2,y_3,y_4\}]$ is a clique. It follows that $g'(y_4)\ge 1.5$ and $\sum_{x\in V_3\cup V_4}g'(x)\ge 3$.
\qed



Since $|V_3\cup V_4|\ge 2$, $\sum_{v\in V_3\cup V_4}g'(v)\ge 2$ by Claim \ref{v34}. When $e(G[V_1])\ge 7$, by inequality (\ref{egprime}), $e(G)\ge 3n-9$. Now we  consider the case $e(G[V_1])= 6$.
 If we can show $\sum_{v\in V_2}g'(v)>0$ or $\sum_{v\in V_3\cup V_4}g'(v)>2$, by  Claim \ref{v34} and  (\ref{egprime}), then $e(G)>3n-10$ and so $e(G)\ge 3n-9$.
 If there exists a vertex $u\in V_2$ such that $|N(u)\cap V_1|\ge 3$, then $g'(u)\ge 1$ and so $\sum_{v\in V_2}g'(v)>0$. If $V_3\ne \emptyset$, then $\sum_{v\in V_3\cup V_4}g'(v)\ge 3$ by Claim \ref{v34}.
Thus we may assume that $|N(v)\cap V_1|=2$ for each $v\in V_2$ and $V_3=\emptyset$.
We choose a vertex $x\in V_2$. Without loss generality, suppose $x\in N(a_1)\cap N(a_2)$. Since $xa_i\notin E$ for each $i\in \{3,4,5\}$, there is a copy of $K_{2,2}$ between $N(a_i)$ and $N(x)$, say $\{a_{i1}, a_{i2}\}\sim \{x_{i1}, x_{i2}\}$.
Note that there is no copy of $K_{1,2}$ in $G[V_1\setminus\{a\}]$. Thus $\{a_{i1},a_{i2}\}\cap V_2\ne \emptyset$ for each $i\in \{3,4,5\}$.
Recall that $|N(v)\cap V_1|=2$ for each $v\in V_2$ and $V_3=\emptyset$. We have  $\{x_{i1},x_{i2}\}\cap (V_2\cup V_4)\ne \emptyset$ for each $i\in \{3,4,5\}$.
By Proposition \ref{k22}(ii), we have $|N(w)\cap V_2|\ge 3$ for each $w\in (\bigcup_{i\in \{3,4,5\}}\{x_{i1},x_{i2}\})\cap (V_2\cup V_4)$. Thus $\sum_{v\in V_3\cup V_4}g'(v)\ge 3$ by Claim \ref{v34} or $\sum_{v\in V_2}g'(v)\ge 0.5$.

Next we  consider $e(G[V_1])=5$. If $\sum_{v\in V\setminus V_1}g'(v)>3$, by (\ref{egprime}), then $e(G)>3n-10$ and so $e(G)\ge 3n-9$. Thus we  prove $\sum_{v\in V\setminus V_1}g'(v)>3$ in the following. 
Recall $\sum_{v\in V_3\cup V_4}g'(v)\ge 2$. If there is a vertex $x\in V_2$ with $|N(x)\cap V_1|\ge 4$, then $g'(x)\ge 2$ and $$\sum_{v\in V\setminus V_1}g'(v)\geq g'(x)+\sum_{v\in V_3\cup V_4}g'(v)\geq 4.$$ If there are two different vertices $x,y\in V_2$ with $|N(x)\cap V_1|=|N(y)\cap V_1|=3$, then $g'(x)\ge 1$,  $g'(y)\ge 1$ and $$\sum_{v\in V\setminus V_1}g'(v)\geq g'(x)+g'(y)+\sum_{v\in V_3\cup V_4}g'(v)\geq 4.$$
Suppose $x\in V_2$ with $|N(x)\cap V_1|=3$, and $|N(v)\cap V_1|=2$ for each $v\in V_2\setminus\{x\}$.  Let $N(x)\cap V_1=\{a_1, a_2, a_3\}$. Since $xa_4\notin E$, there is a copy of $K_{2,2}$ between $N(x)$ and $N(a_4)$, say $\{x_{11},x_{12}\}\sim \{a_{41},a_{42}\}$. If $V_3\ne \emptyset$, by Claim \ref{v34}, then $\sum_{v\in V_3\cup V_4}g'(v)\ge 3$. Thus $$\sum_{v\in V\setminus V_1}g'(v)\ge g'(x)+\sum_{v\in V_3\cup V_4}g'(v)\geq4.$$
So we may assume that $V_3=\emptyset$. Then $\{a_{41},a_{42}\}\subseteq V_2$ and $\{x_{11},x_{12}\}\cap (V_2\cup V_4)\ne \emptyset$.
Let $w\in \{x_{11},x_{12}\}\cap (V_2\cup V_4)$. By Proposition \ref{k22}(ii), $|N(w)\cap V_2|\ge3$. Thus $g'(w)\ge 0.5$. If $w\in V_2$,
 then $$\sum_{v\in V\setminus V_1}g'(v)\ge g'(x)+g'(w)+ \sum_{v\in V_3\cup V_4}g'(v)\geq 3.5.$$ If $w\in V_4$, by Claim \ref{v34}, then $\sum_{v\in V\setminus V_1}g'(v)\ge 4$.

Suppose $|N(v)\cap V_1|=2$ for each $v\in V_2$.
Since $|V_3\cup V_4|\le 6$ and $n\geq 20$, $|V_2|\ge 8$. Recall the definition of  $g$-function, for each $v\in V_2$, we have $g(v)\ge 0.25$ and if $g(v)>0.25$, then $g(v)\ge 0.5$.
We see there exists a vertex $x\in V_2$ such that $g(x)=0.25$, otherwise, $g(v)\ge 0.5$ for each $v\in V_2$ and so  $\sum_{v\in V_2}g(v)\ge 0.5|V_2|\ge 4$.
By  (\ref{eg}), $e(G)\ge 5+3(n-6)+4=3n-9$.
We choose such a vertex $x\in V_2$ such that $g(x)=0.25$. Then $d(x)=5$ and let $N(x)=\{a_1,a_2,x_{11},x_{12},z\}$, where $\{a_1, a_2\}\subseteq V_1$, $\{x_{11},x_{12}\}\subseteq V_2$ and $z\in V_4$. Note that $xa_j\notin E$ for each $j\in \{3,4,5\}$. By Proposition \ref{k22}(i), there is a copy of $K_{2,2}$ between $N(x)$ and $N(a_j)$, say $\{x_{j1},x_{j2}\}\sim \{a_{j1},a_{j2}\}$. We see $\{a_{j1},a_{j2}\}\subseteq V_2\cup V_3$ for each $j\in \{3,4,5\}$. Since $|N(v)\cap V_1|=2$ for each $v\in V_2$, $\{x_{j1},x_{j2}\}\nsubseteq V_1$ for each $j\in \{3,4,5\}$. Otherwise, $\{x_{j1},x_{j2},a_j\}\subseteq N(a_{j1})\cap V_1$, a contradiction.

Suppose $V_3\neq \emptyset$. Then we  have $\sum_{v\in V_3\cup V_4}g'(v)\ge 3$ by Claim \ref{v34}.
We have $|V_3|\le 3$, otherwise $\sum_{v\in V_3}g'(v)\ge 4$ and we are done.  Note that $\{a_{j1},a_{j2}\} \subseteq V_2\cup V_3$ for any $j\in \{3,4,5\}$. When $|\bigcup_{j\in \{3,4,5\}}\{a_{j1},a_{j2}\}|\leq 5$, we may assume $a_{31}=a_{41}$.  When $|\bigcup_{j\in \{3,4,5\}}\{a_{j1},a_{j2}\}|=6$, we have $|\bigcup_{j\in \{3,4,5\}}\{a_{j1},a_{j2}\}\cap V_2|\geq 3$ because $|V_3|\le 3$, so we may assume $\{a_{31}, a_{41}\}\subseteq \bigcup_{j\in \{3,4,5\}}\{a_{j1},a_{j2}\}\cap V_2$.
In two cases, we have $\{a_{31}, a_{41}\}\subseteq V_2$.
Let $k\in\{3,4\}$. If $\{x_{k1},x_{k2}\}\cap V_2\ne \emptyset$, let $w\in \{x_{k1},x_{k2}\}\cap V_2$, then
$\{x,a_{k1}\}\subseteq N(w)\cap V_2$, Proposition \ref{k22}(ii) implies that $|N(w)\cap V_2|\ge 3$ and so $g'(w)\ge 0.5$. Thus $$\sum_{v\in V\setminus V_1}g'(v)\ge g'(w)+ \sum_{v\in V_3\cup V_4}g'(v)\geq3.5.$$
So we assume  $\{x_{k1},x_{k2}\}\cap V_2=\emptyset$.
Note that $\{x_{k1},x_{k2}\}\nsubseteq V_1$.
Since $d(x)=5$,
 $\{x_{k1},x_{k2}\}=\{a_1,z\}$ or $\{a_2,z\}$, and so $N(a_{k1})\cap N(a_{k2})\cap V_1=\{a_k,a_{\ell_k}\}$ for some $\ell_k\in [2]$.  Then $\{x,a_{31},a_{32},a_{41},a_{42}\} \subseteq N(z)\cap V_2$. Note that $|N(v)\cap V_1|=2$ for any $v\in V_2$. Since
  $x\in V_{12}$, $\{a_{31},a_{32}\}\subseteq V_{1\ell_3}$ and $\{a_{41},a_{42}\}\subseteq V_{1\ell_4}$,  $|\{x,a_{31},a_{32},a_{41},a_{42}\}|= 5$, which follows that
 $g'(z)\ge 2$.  Note that $V_3\ne \emptyset$ and  $g'(y)\ge 1$ for each $y\in V_3$ and $g'(v)\ge 0.5$ for each $v\in V_4$. Recall $z\in V_4$ and $|V_3\cup V_4|\ge 2$.
 So $\sum_{v\in V_3\cup V_4}g'(v)>3$ when $|V_3\cup V_4|\geq3$. If $|V_3\cup V_4|=2$, then $g'(y)>1$ for $y\in V_3$ because $d(y)\ge 5$. Therefore $\sum_{v\in V_3\cup V_4}g'(v)>3$.

It remains to consider $V_3=\emptyset$. Then $\{a_{j1},a_{j2}\}\subseteq V_2$ for any $j\in \{3,4,5\}$. When $\{x_{j1},x_{j2}\}\cap V_2=\emptyset$ for any $j\in \{3,4,5\}$, then $\{x_{j1},x_{j2}\}=\{a_1,z\}$ or $\{a_2,z\}$. Note that $|N(v)\cap V_1|=2$ for each $v\in V_2$. Since $\{a_{j1},a_{j2}\}\subseteq V_{j\ell_j}$ for  $\ell_j\in [2]$, $|(\bigcup_{j\in \{3,4,5\}}\{a_{j1},a_{j2}\})\cup \{x\}|=7$. Thus $|N(z)\cap V_2|\ge 7$, which implies that $\sum_{v\in V_2}g'(v)\ge 4$.
When there exists $j\in \{3,4,5\}$ such that $\{x_{j1},x_{j2}\}\cap V_2\ne \emptyset$, then $g'(w)\ge 0.5$ for $w\in   \{x_{j1},x_{j2}\}\cap V_2$ because $|N(w)\cap V_2|\ge 3$ by Proposition \ref{k22}(ii). In this case, we have $z\notin \{x_{j1},x_{j2}\}\cap V_4$. Otherwise,  Proposition \ref{k22}(ii) implies $|N(z)\cap V_2|\ge 3$. By
 Claim \ref{v34}, $$\sum_{v\in V\setminus V_1}g'(v)\ge g'(w)+\sum_{v\in V_3\cup V_4}g'(v)\geq 3.5.$$ Thus we are done.
If $|N(x_{11})\cap V_2|+|N(x_{12})\cap V_2|\ge 7$, then we have $$g'(x_{11})+g'(x_{12})= e(G[\{x_{11},x_{12}\},V_1])+0.5(e(G[\{x_{11}\},V_2])+e(G[\{x_{12}\},V_2]))-6\ge 4+3.5-6=1.5.$$ Thus $\sum_{v\in V\setminus V_1}g'(v)\ge 3.5$, and we are done. So it suffices to prove $|N(x_{11})\cap V_2|+|N(x_{12})\cap V_2|\ge 7$ in the following.
Since $z\notin \{x_{j1},x_{j2}\}$, we have $\{x_{j1},x_{j2}\}\cap V_2\ne \emptyset$ for any $j\in \{3,4,5\}$. 
Recall $N(x)=\{a_1, a_2, x_{11}, x_{12}, z\}$ and $x\in N(x_{11})\cap N(x_{12})\cap V_{12}$. Then $\{a_{31},a_{32},a_{41},a_{42},a_{51},a_{52}\}\subseteq N(x_{11})\cup N(x_{12})$.
If $|\{a_{31},a_{32},a_{41},a_{42},a_{51},a_{52}\}|\ge 5$, then
$$|N(x_{11})\cap V_2|+|N(x_{12})\cap V_2|=|(N(x_{11})\cup N(x_{12}))\cap V_2|+|(N(x_{11})\cap N(x_{12}))\cap V_2|
\ge 7.$$
Suppose that $|\{a_{31},a_{32},a_{41},a_{42},a_{51},a_{52}\}|\le 4$.
Note that $|N(x)\cap V_1|=2$ for each $x\in V_2$. We obtain  $|\{a_{31},a_{32},a_{41},a_{42},a_{51},a_{52}\}|\ge 3$.
When $\{x_{31},x_{32}\}\cap V_1\ne \emptyset$, say $a_\ell\in \{x_{31},x_{32}\}$ for some $\ell\in [2]$, then $\{a_{31},a_{32}\} \subseteq V_{3\ell}$
 and  $\{a_{31},a_{32}\}\cap \{a_{k1},a_{k2}\}=\emptyset$ for each $k\in \{4,5\}$.
Since $|\{a_{31},a_{32},a_{41},a_{42},a_{51},a_{52}\}|\le 4$, we have $\{x_{k1},x_{k2}\}=\{x_{11},x_{12}\}$ for   each $k\in \{4,5\}$, that is $|N(x_{11})\cap N(x_{12}) \cap \bigcup_{j\in \{3,4,5\}}\{a_{j1},a_{j2}\}|\ge 2$.
Thus
\begin{align*}
|N(x_{11})\cap V_2|+|N(x_{12})\cap V_2|=&\notag|(N(x_{11})\cup N(x_{12}))\cap V_2|+|(N(x_{11})\cap N(x_{12}))\cap V_2|\\\notag
\ge&|\bigcup\nolimits_{j\in \{3,4,5\}}\{a_{j1},a_{j2}\}\cup \{x\}|+3 \geq 7.
\end{align*}
When $\{x_{j1},x_{j2}\}\cap V_1= \emptyset$ for each $j\in\{3,4,5\}$, then $\{x_{j1},x_{j2}\}=\{x_{11},x_{12}\}$ for each $j\in \{3,4,5\}$ and $\bigcup_{j\in \{3,4,5\}}\{a_{j1},a_{j2}\}\subseteq N(x_{11})\cap N(x_{12})$.
By $|\bigcup_{j\in \{3,4,5\}}\{a_{i1},a_{i2}\}|\ge 3$ and $x\notin\bigcup_{j\in \{3,4,5\}}\{a_{i1},a_{i2}\}$, we have $|N(x_{11})\cap V_2|+|N(x_{12})\cap V_2|\geq 8$.


As a result, we have $e(G)\ge 3n-9$ for  $n\ge 9$ in each case and so $sat_5(n,K_{3,3})\ge 3n-9$. 
\qed
  
This completes the proof of Theorem \ref{main23}.
\section{Conclusion}

Based on above results, we  make  the following conjecture, which  gives an exact value for $sat(n,K_{3,3})$.

\begin{conj}\label{k33conj}
	For $n\ge 9$, $sat(n,K_{3,3})= 3n-9$.
\end{conj}

By Theorem \ref{main1}, $sat(n,K_{3,3})\le 3n-9$ for $n\ge 9$. To confirm Conjecture \ref{k33conj}, it suffices to  prove $sat(n,K_{3,3})\ge 3n-9$ for $n\ge 9$.
Let $G$ be a $K_{3,3}$-saturated graph with $n$ vertices and $n\geq9$. Proposition \ref{k22}(i) implies $\delta(G)\ge 2$. If $\delta(G)\ge 6$, then $e(G)\ge 3n\geq 3n-9$.
Thus we only need to consider $2\le \delta(G)\le 5$. We have proved $sat_{\delta}(n,K_{3,3})\ge 3n-9$ when $\delta\in\{2,5\}$.
Actually, for $\delta\in \{3,4\}$, we  can also apply the method in this paper, but it is more complex and there are quite a few  cases to consider.

\bigskip

{\bf Acknowledgments.} Shi and Zhang are partially supported  by the National Natural Science Foundation of China (Nos. 11922112, 12161141006), Natural Science Foundation of Tianjin (Nos. 20JCZDJC00840, 20JCJQJC00090).  Lei was partially supported by the National Natural Science Foundation of China (No. 12001296), Natural Science Foundation of Tianjin (No. 21JCQNJC00060).

\end{document}